\documentclass[11pt,a4paper]{article}

\usepackage{amsmath, amsfonts, amssymb}
\usepackage{color}
\usepackage{graphicx}

\def\be#1\ee{\begin{equation}#1\end{equation}}

\newtheorem{thm}{Theorem}

\newtheorem{prop}[thm]{Proposition}

\newtheorem{rem}[thm]{Remark}

\def\P{{\mathbb{P}}}
\def\R{\mathbb{R}}

\def\AA{{\mathcal A}}

\def\LL{{\mathcal{L}}}
\def\MM{{\mathcal{M}}}

\def\OO{{\Theta}}
\def\OOO{{\mathcal{O}}}
\def\PP{{\mathcal P}}

\let\BFseries\bfseries\def\bfseries{\BFseries\mathversion{bold}} 

\title{On Brownian Exit Times \\
  from Some Non-Convex Domains }
\author{
M.\,Lifshits\footnote{St.Petersburg State University, Russian Federation, 199034,
St.Petersburg, Universitetskaya emb. 7-9; E-mail: m.lifshits@spbu.ru
}
\ and
A.\,Nazarov\footnote{ St.Petersburg Department of Steklov Mathematical Institute of
Russian Academy of Science, Russian Federation, 191023, St.Petersburg, Fontanka emb. 27
and St.Petersburg State University, Russian Federation, 191023, St.Petersburg,
Universitetskaya emb. 7-9; E-mail: nazarov@pdmi.ras.ru
}
}
\bigskip

\begin{document}

\maketitle

\begin{abstract}
    We consider the  tail probabilities for Brownian exit time from a class
    of {\it perturbed multi-strips} in Euclidean plane. Under some assumptions
    we prove that the long stays in a perturbed multi-strip are more likely than
    those in a strip of the same width. This effect is related to the existence
    of the trapped modes in waveguides.
\end{abstract}

\medskip

 {\bf AMS classification codes:}  Primary: 60J65. Secondary: 35K05, 35J25.

\section{Introduction and results}
\medskip

Let $W(t)=\big(W_j(t)\big)\big|_{j=1}^d$, $t\ge 0$, be a standard $d$-dimensional
Brownian motion with independent coordinates, $d\ge 1$, and let $W^x(t)=x+W(t)$
denote a Brownian motion starting from a point $x\in \R^d$.

We begin from two natural small deviation problems: find the asymptotic behavior of
\begin{eqnarray*}
   \P( \| \min_{1\le j\le d} |W_j|\| \le r) , \qquad r \to 0,
\\
   \P( \| \min_{1\le j\le d} W_j\| \le r  ), \qquad r \to 0.
\end{eqnarray*}
Here $\|\cdot\|$ stands for sup-norm, i.e. $\|X\|=\max_{0 \le t \le 1} |X(t)|$.

\begin{figure*}[ht]
\label{f:xxx}
\centering
\includegraphics[width=8cm]{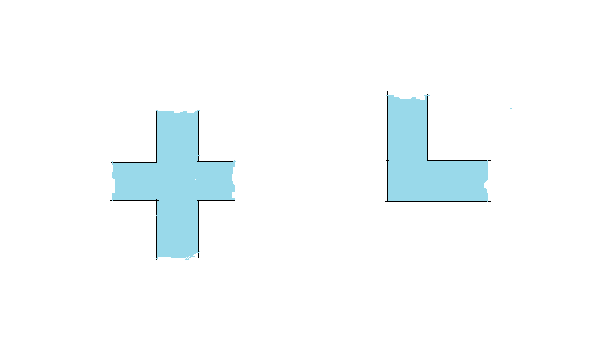}
\caption{ Two-dimensional cross and two-dimensional corner}
\label{figureMAL1}
\end{figure*}

This problem immediately reduces to the so-called exit time asymptotics. Denote by
$\tau_{x,\OO}$ the exit time for $W^x$ from a domain $\OO\subset\R^d$ containing
a point $x$, and put $\tau_{\OO}=\tau_{0,\OO}$. Then, by self-similarity of $W$,
we have
\begin{eqnarray*}
   \P(  \| \min_{1\le j\le d} |W_j|\| \le r)
    &=& \P\left\{ \tau_{\OO_1} \ge r^{-2}\right\},
\\
   \P( \| \min_{1\le j\le d} W_j\| \le r )
   &=& \P\left\{ \tau_{\OO_2} \ge r^{-2}\right\},
\end{eqnarray*}
where
\[
  \OO_1:= \{ x\in\R^d: \min_{1\le j\le d} |x_j| < 1 \}
\]
is a cross, and
\[
  \OO_2:= \{ x\in\R^d: |\min_{1\le j\le d} x_j| < 1 \}
\]
is a corner, see Fig.~\ref{f:xxx} in dimension $2$.

It is well known that the function
\be \label{uxt}
  u(x,t):= \P\{W^x(s)\in \OO, 0\le s < t\} = \P\{\tau_{x,\OO}\ge t\}
\ee
solves the initial-boundary value problem
\begin{equation} \label{heat}
   \LL u\equiv\partial_tu-\frac 12\Delta u=0
   \quad\mbox{in}\quad \OO\times \R_+; \qquad  u\big|_{x\in\partial\OO}=0;
   \qquad u\big|_{t=0}\equiv 1.
\end{equation}

The solution of \eqref{heat} for $d=1$, $\Theta=(-1,1)$,
\begin{equation} \label{heat1}
    u(x,t) =\frac{4}{\pi}\, \sum_{k=0}^\infty
    \frac{(-1)^k}{2k+1}\, \cos\left(\pi(k+1/2)x\right) \exp({-(k+1/2)^2\pi^2t/2}),
\end{equation}
known long ago, was used in the one-dimensional exit probability problem
by Chung, cf. \cite{Chung} or \cite[v.2, Ch. X]{Feller}.
Formula \eqref{heat1} yields the asymptotics
\begin{equation*} \label{tau1}
  \P\{\tau_{x,\OO}\ge t\} = \frac{4}{\pi}\,\cos(\pi x/2)  \exp(-\pi^2 t/8)+ O\big(\exp(-9\pi^2 t/8)\big) ,
    \quad \textrm{as } t\to \infty.
\end{equation*}
Since any one-dimensional projection of a multivariate Brownian motion is
a one-dimensional Brownian motion, we obtain the same results when $d>1$ and
\[
  \OO=\Pi:=\{x\in\R^d: |x_d| < 1\}
\]
is a layer (strip for $d=2$). For the exit times from the semi-strip in $\R^2$
\[
   \PP := \{x\in\R^2: x_1>0, \ |x_2|< 1\}
\]
we have essentially the same tail asymptotics because for
$x:=(x_1,x_2)\in \PP$ one has
\begin{eqnarray*}
  \P(\tau_{x,\PP}\geq t)&=& \P(\tau_{x_1,\R_+}\geq t)\,
  \P(\tau_{x_2,(-1,1)}\geq t)
\\
  &\sim& \sqrt{\frac{2}{\pi t}}\,x_1 \cdot\,
  \frac{4}{\pi}\,\cos(\pi x_2/2)  \exp({-\pi^2 t/8}),
  \quad \textrm{as } t\to \infty.
\end{eqnarray*}
One may also cite the works \cite{banuelos-deblassie-smits,banuelos-smits,WLi,LifZS}
investigating the tails of $\tau_{\OO}$ for conic, parabolic,
and some other unbounded convex domains of this kind.
\medskip

Notice that domains $\OO_1$ and $\OO_2$ for $d=2$ look like a bunch of
semi-strips when one observes them far from the origin. More formally,
we say that a domain $\OO$ in $\R^2$  is a {\it perturbed multi-strip}
if outside of some compact $\OO$ coincides with a finite union
of non-intersecting domains $\PP_j$ congruent to $\PP$.

From the first glance, it might seem that the decay rate of the tail probability
of the exit time from a perturbed multi-strip should be the same as for $\PP$.
Somewhat surprisingly it is not the case in general. For many domains of this type
(in particular, for two-dimensional cross and corner) these tail probabilities
decay more slowly.

The reason of this phenomenon is as follows. It is well known that the spectrum
of the Dirichlet Laplacian in the strip is absolutely continuous:
$$
   \sigma(-\Delta,\Pi)=\sigma_{ac}(-\Delta,\Pi)=[\Lambda,\infty),
$$
where $\Lambda=\pi^2/4$ is the first eigenvalue of one-dimensional Dirichlet
Laplacian in the cross-section $(-1,1)$.

The same is true in the semi-strip $\PP$. This gives us the {\it essential
spectrum} of the Dirichlet Laplacian in any perturbed multi-strip $\OO$:
$$
   \sigma_{ess}(-\Delta,\OO)=\sigma_{ac}(-\Delta,\OO)=[\Lambda,\infty),
$$
see, e.g., \cite[Proposition 1.1.1]{EK} for a more general case.

However, this does not exclude the existence of {\it discrete spectrum} of
the Dirichlet Laplacian in $\OO$ below the threshold $\Lambda$.
This is indeed the case for two-dimensional cross and corner, see \cite{SRW},
\cite[Proposition 1.5.2]{EK} and \cite{He}, \cite[Proposition 1.2.2]{EK},
respectively.
We show below that this discrete spectrum is responsible for the slower decay
of the tail probability of the exit time.

Notice that the phenomenon of emerging additional discrete spectrum in
perturbed strips is deeply studied in theory of {\it waveguides} where it
generates so-called {\it trapped modes}, see, e.g., \cite{EK} and references
therein.
\medskip

Thus, in what follows we assume that the perturbed multi-strip $\OO$ satisfies\medskip

{\bf Assumption A}. The spectrum of the Dirichlet Laplacian in $\OO$ consists
of the ray $[\Lambda,\infty)$
and a non-empty discrete spectrum in the interval $(0,\Lambda)$.\medskip

We denote by $\lambda_0$ the smallest eigenvalue and by $\lambda_1$ the infimum
of other spectrum (note that $\lambda_1$ either coincides with $\Lambda$ or is
another eigenvalue).

\begin{prop}\label{prop1}
Let Assumption A be satisfied for a perturbed multi-strip $\OO$.
Then the eigenvalue $\lambda_0$ is simple, and corresponding eigenfunction
$v_0$ does not vanish on $\OO$ and decays exponentially as $x\to\infty$. Namely,
if we direct the axis $y_1$ along the axis of semi-strip $\PP_j$, then for $y\in\PP_j$
\begin{equation}  \label{v0}
     v_0(y)= C_j\exp ( -\sqrt{\Lambda-\lambda_0} \,y_1)\cos(\sqrt{\Lambda}\,y_2)\cdot(1+o(1)),
     \quad |y|\to\infty.
\end{equation}
\end{prop}

This fact is well known. We give a sketch of the proof for the reader's
convenience.

The smallest eigenvalue of the Dirichlet Laplacian can be found as minimum
of the Rayleigh quotient
$$
   Q(v)= \frac {\int\limits_{\OO}|\nabla v|^2 } {\int\limits_{\OO}|v|^2}
$$
over the Sobolev space $\stackrel {o}{W}\!\!\vphantom{W}^1_2(\OO)$, and
corresponding eigenfunction minimizes this Rayleigh quotient.
Since $Q(|v|)=Q(v)$, for arbitrary eigenfunction $\widehat v$ corresponding
to $\lambda_0$, $|\widehat v|$ is also an eigenfunction.
Therefore, $|\widehat v|$ is superharmonic, and by the maximum principle it cannot
vanish inside $\OO$.

Next, if $\lambda_0$ were a multiple eigenvalue, we could choose a pair of
eigenfunctions orthogonal in $L_2(\OO)$. Then at least one of them
should change sign in $\OO$ that is impossible by the previous argument.

Finally, the relation \eqref{v0} can be easily obtained for every semi-strip
$\PP_j$ using the separation of variables in the equation for $v_0$.\hfill$\square$\medskip

Without loss of generality we assume $v_0$ normalized in $L_2(\OO)$.\medskip

Our main result is as follows.

\begin{thm} \label{t1} Let a perturbed multi-strip $\OO$ satisfy Assumption A.
Then the tail probability $\eqref{uxt}$ admits a representation
$$
  u(x,t)=\AA \, v_0(x)\exp(-\frac {\lambda_0}2 t)
   +\widehat u(x,t),\quad |\widehat u(x,t)|\le C(t+1)\exp(-\frac {\lambda_1}2 t),
$$
where $\AA = \int\limits_{\OO} v_0$  while $\lambda_0<\lambda_1\le \Lambda$
are introduced after Assumption A.
\end{thm}

\begin{rem}
 We conjecture that in fact the remainder term admits a better estimate $|\widehat u(x,t)|\le C\exp(-\frac {\lambda_1}2 t)$.
\end{rem}

We give the proof of Theorem \ref{t1} in Section 2.
In Section 3 we discuss this result and give some generalizations.

\section{Proof of Theorem \ref{t1}}

In what follows we use a standard notation
$\big(v,w\big)_{\OO}:=\int\limits_{\OO}vw$.
Various constants depending only on $\OO$ are denoted by $C$.

A typical example of a  perturbed multi-strip $\OO$ is a connected union
of the strip $\Pi\subset \R^2$ and a bounded domain
$\OOO\subset \{x\in\R^2: |x_1|< R\}$. If $\OO\ne\Pi$ then
Assumption A is satisfied for such $\OO$, see e.g. \cite[Theorem 1.4]{EK}.
We prove the statement for this $\OO$, the proof for a general case follows
the same line.

Let $U(x,t)$ be the solution of the corresponding problem in the strip $\Pi$:
\begin{equation*}
    \LL U=0 \quad\mbox{in}
    \quad \Pi\times \R_+; \qquad  U\big|_{x\in\partial\Pi}=0;
    \qquad U\big|_{t=0}\equiv 1.
\end{equation*}
Obviously, $U$ does not depend on $x_1$, cf. \eqref{heat1}, and
$0<U(\cdot,t)<C\exp(-\frac {\Lambda}2 t)$.

Now we split $u$ into two parts: $u=u_1+u_2$, where
$u_1(x,t)=U(x_2,t)\chi(x_1)$ and $\chi$ is a smooth cutoff function:
$$
   \chi(x_1)=0
    \ \ \mbox{for} \ \ |x_1|<R;
    \quad \chi(x_1)=1\ \ \mbox{for}
    \ \ |x_1|>R+1.
$$
Notice that $\varphi=u_2(\cdot,0)$ has compact support, and
the support of the function
$$
   f:=\LL u_2=-\LL u_1=\frac 12 UD^2_{x_1}\chi
$$
is compact w.r.t. $x$. Therefore, we can apply the spectral
decomposition for the operator exponent:
\begin{equation}\label{decomp}
    \aligned
&    u_2(\cdot,t) = \exp(\frac 12\Delta t)\varphi(\cdot)
     + \int\limits_0^t\exp(\frac 12\Delta (t-s))f(\cdot,s)\,ds
    \\
&    =\! \int\limits_{[\lambda_0,\infty)}\!\! \exp(-\frac {\lambda}2 t)\,dE(\lambda)\varphi(\cdot)
      +\int\limits_0^t\int\limits_{[\lambda_0,\infty)}\!\!\exp(-\frac {\lambda}2 (t-s))
      \,dE(\lambda)f(\cdot,s)\,ds,
      \endaligned
\end{equation}
where $E(\lambda)$ is the projector-valued spectral measure generated
by the (self-adjoint) Dirichlet-Laplacian operator in $L_2(\OO)$,
see \cite[Ch. 6]{BS}.

Since $\lambda_0$ is a simple eigenvalue, we can rewrite \eqref{decomp}
as $u_2=u_{2,0}+u_{2,1}$, where
\begin{equation*}
    u_{2,0}(\cdot,t) = \Big(\exp(-\frac {\lambda_0}2 t)
    \big(\varphi,v_0\big)_{\OO}
    +\int\limits_0^t\exp(-\frac {\lambda_0}2 (t-s))
    \big(f(\cdot,s),v_0\big)_{\OO}\,ds\Big) v_0(\cdot);
\end{equation*}
\begin{eqnarray*}
    u_{2,1}(\cdot,t) &=& u_{2,11}(\cdot,t)+u_{2,12}(\cdot,t)
    :=\int\limits_{[\lambda_1,\infty)}\exp(-\frac {\lambda}2 t)\,dE(\lambda)
    \varphi(\cdot)
\\
    &+&\int\limits_0^t\int\limits_{[\lambda_1,\infty)}\!
      \exp(-\frac {\lambda}2 (t-s))\,dE(\lambda)f(\cdot,s)\,ds.
\end{eqnarray*}

%

First, we estimate $L_2$-norm of $u_{2,1}$. Estimate
of the first term is evident:
$$
	\|u_{2,11}(\cdot,t)\|^2_{L_2(\OO)}
    \le\exp(-\lambda_1 t)\,\|\varphi\|^2_{L_2(\OO)}.
$$
Since $|f(\cdot,s)|<C\exp(-\frac {\Lambda}2 s)$ and $f(\cdot,s)$ is
compactly supported uniformly w.r.t. $s$, the estimate for
$u_{2,12}$ follows from $\lambda_1\le\Lambda$:
$$
	\|u_{2,12}(\cdot,t)\|^2_{L_2(\OO)}\le C\exp(-\lambda_1 t)\,
	\Big(\int\limits_0^t\exp(\frac {\lambda_1-\Lambda}2\, s)\,ds\Big)^2
    \le Ct^2\exp(-\lambda_1 t).
$$
Thus, we have $\|u_{2,1}(\cdot,t)\|_{L_2(\OO)}\le C(t+1)\exp(-\frac {\lambda_1}2 t)$.
Therefore, standard parabolic estimates give
$|u_{2,1}(\cdot,t)|\le C(t+1)\exp(-\frac {\lambda_1}2 t)$.\medskip

Next, we transform the second term in $u_{2,0}$ using integration
by parts:
\begin{eqnarray*}
\nonumber
   \big(f(\cdot,s),v_0\big)_{\OO}
   &=& \frac 12\big(U(\cdot,s)D^2_{x_1}\chi,v_0\big)_{\OO}
   =\frac 12\big(U(\cdot,s)\chi,D^2_{x_1}v_0\big)_{\OO}
\\
   &=& -\frac {\lambda_0}2\big(U(\cdot,s)\chi,v_0\big)_{\OO}
   -\frac 12\big(U(\cdot,s)\chi,D^2_{x_2}v_0\big)_{\OO}
\\
   &=& -\frac {\lambda_0}2\big(U(\cdot,s)\chi,v_0\big)_{\OO}
   -\frac 12\big(D^2_{x_2}U(\cdot,s)\chi,v_0\big)_{\OO}
\\
   &=& -\frac {\lambda_0}2\big(U(\cdot,s)\chi,v_0\big)_{\OO}
   -\big(\partial_sU(\cdot,s)\chi,v_0\big)_{\OO}
\end{eqnarray*}
(notice that all integrals converge due to the exponential
decay of $v_0$);
\begin{eqnarray*}
  &&\int\limits_0^t\exp(-\frac {\lambda_0}2 (t-s))
     \big(f(\cdot,s),v_0\big)_{\OO}\,ds
\\
  &=& -\exp(-\frac {\lambda_0}2 t) \int\limits_0^t
    \exp(\frac {\lambda_0}2 s)
     \big(\frac {\lambda_0}2U(\cdot,s)\chi
     +\partial_sU(\cdot,s)\chi,v_0\big)_{\OO}\,ds
\\
  &=& -\exp(-\frac {\lambda_0}2 t)\cdot
  \big(\exp(\frac {\lambda_0}2 s)U(\cdot,s)
  \chi,v_0\big)_{\OO}\Big|_0^t
  \\
  &=& \exp(-\frac {\lambda_0}2 t)\cdot
  \big(U(\cdot,0)\chi,v_0\big)_{\OO}
  -\big(U(\cdot,t)\chi,v_0\big)_{\OO}.
\end{eqnarray*}

We sum up all terms and take into account that
$\varphi+U(\cdot,0)\chi\equiv1$. This gives
$$
   u(x,t)=\big(1,v_0\big)_{\OO}\,v_0(x)
   \exp(-\frac {\lambda_0}2 t)
   +
   O\Big((t+1)\exp(-\frac {\lambda_1}2 t)\Big),
$$
as required.\hfill$\square$

\section{Discussion and generalizations}

We begin with some comments.

\begin{enumerate}
\item The relation $\lambda_0 < \Lambda$ shows that the long stays
in a perturbed multi-strip under Assumption A are more likely than
those in a strip of the same width.
\item By Proposition \ref{prop1}, the main term of the tail
probability \eqref{uxt} depends on the initial point $x$ exponentially.
This is a new phenomenon that just does not exist in dimension one.
If we start from a remote $x$, the optimal strategy to stay in $\OO$
for a longer time is to run towards the origin and then stay somewhere
near it. The relation between the running time and staying time is a
subject of optimization, which results in \eqref{v0}. Notice that
there is no "running" phase for long stays in a strip.
\item For the case of cross we can give a simple geometric explanation
of the inequality $\lambda_0 < \Lambda$. Consider the square
$\widetilde\OO_1:=\{x\in\R^2:\, |x_1| + |x_2|< 2 \}$. Using the version
of Chung's formula \eqref{heat1} for squares, we easily obtain for
$\widetilde\OO_1$ the same rate of decay as for the strip. On the other
hand, since $\widetilde\OO_1\subsetneq \OO_1$, the inequality between
the main exponents of the two sets is strict. For the case of corner we
do not know an elementary proof.

Numerical calculations give for the cross in
$\R^2$ $\lambda_0\approx 0.66\Lambda$, see, e.g.,
\cite[Proposition 1.5.2]{EK},
and for the corner  in $\R^2$ $\lambda_0\approx 0.929\Lambda$, see,
e.g., \cite{ESS}, \cite[Proposition 1.2.3]{EK},
For both domains it is also known that the eigenvalue under threshold
is unique.
\end{enumerate}

Next, we notice that the proof of Theorem \ref{t1} runs with minor
changes for two classes of domains in $\R^d$, $d\ge3$:

\begin{enumerate}
 \item A connected union of the layer $\Pi\subset \R^d$ and a bounded
 domain $\OOO\subset \{x\in\R^d: |x|< R\}$.
 In this case $\Lambda=\pi^2/4$, and Assumption A for $\OO\ne\Pi$
 follows from \cite[Theorem 4.5]{EK}.

\item A {\it perturbed multi-tube} that is a domain in $\R^d$ such
 that outside of a compact it coincides with a finite union
 of non-intersecting domains congruent to a semitube
$$
   \widehat\PP:= \{x\in\R^d: x_1>0, \ (x_2,\dots,x_d)\in\MM\},
$$
 where $\MM$ is a bounded domain in $\R^{d-1}$ with ``not very bad''
(say, Lipschitz) boundary. In this case
 $\Lambda$ is the first eigenvalue of $(d-1)$-dimensional Dirichlet
 Laplacian in the cross-section $\MM$.

The Assumption A should be verified separately. For instance, in the
paper \cite{BMN} the existence and uniqueness of the eigenvalue below
the threshold was established for the union of two circular unit
 cylinders whose axes intersect at the right angle.
\end{enumerate}

As for original small deviation problems in dimensions $d\ge3$, they
require more thorough study of the corresponding spectral structure
and will be investigated in a forthcoming paper.

\bigskip

{\bf Acknowledgement.} This research was supported by Russian Foundation of
Basic Research grant 16-01-00258 and by the co-ordinated
grants of DFG (GO420/6-1) and St.\ Petersburg State University (6.65.37.2017).

We are grateful to E.~Hashorva for pointing us this problem and to F.~Bakharev for valuable discussion and references.

\end{document}